\documentstyle{amsppt}
\pagewidth{6in}\vsize8.5in\parindent=6mm\parskip=3pt\baselineskip=14pt
\tolerance=10000\hbadness=500
\NoRunningHeads
\loadbold
\topmatter
\title
Mean square discrepancy bounds for the number of lattice points in
large convex bodies
\endtitle
\author
Alexander Iosevich \ \ \  Eric Sawyer \ \ \
Andreas Seeger
\endauthor
\thanks Research supported in part by  NSF grants.  \endthanks
\address
A. Iosevich,
Mathematics Department, University of Missouri, Columbia, MO 65211, USA
\endaddress
\email iosevich\@math.missouri.edu \endemail
\address
E. Sawyer,
Department of Mathematics and Statistics, 1280 Main Street West,
 Hamilton, Ontario  L8S 4K1, Canada
\endaddress
\email
sawyer\@mcmaster.ca \endemail
\address
A. Seeger, Department of
 Mathematics,
University of Wisconsin-Madison,
Madison, WI 53706, USA
\endaddress
\email seeger\@math.wisc.edu \endemail
\dedicatory Dedicated to the memory of Thomas Wolff\enddedicatory
\endtopmatter

\def\leaderfill{\leaders\hbox to 1em{\hss.\hss}\hfill}

\def\volume{\text{vol}}
\def\R{{\Bbb R}}

\define\dist{{\text{\rm dist}}}

\define\supp{{\text{\rm supp }}}

\define\inn#1#2{\langle#1,#2\rangle}

\define\lcontr{\rfloor}
\define\lco#1#2{{#1}\lcontr{#2}}
\define\lcoi#1#2{\imath({#1}){#2}}
\define\rco#1#2{{#1}\rcontr{#2}}
\redefine\exp{{\text{\rm exp}}}
\define\bin#1#2{{\pmatrix {#1}\\{#2}\endpmatrix}}

\define\card{\text{\rm card}}
\define\lc{\lesssim}
\define\gc{\gtrsim}
             \define\Ga{\Gamma}
\define\eps{\varepsilon}
\define\veps{\epsilon}

              \define\Si{\Sigma}

              \define\Om{\Omega}
\define\ka{\kappa}

\define\bbR{{\Bbb R}}
\define\bbZ{{\Bbb Z}}

\define\cC{{\Cal C}}

\define\cG{{\Cal G}}

\define\cN{{\Cal N}}

\define\cU{{\Cal U}}

\def\cZ{{\Cal Z}}

\define\ic{\imath}

\define\fS{{\frak S}}

\define\fm{\frak m}

\head {\bf 1. Introduction}\endhead
 Let $\Omega$ be a convex domain in $\bbR^d$ containing the origin in its interior. 
We mostly assume that $\Omega$ has  smooth boundary and
 that the Gaussian curvature of the boundary vanishes nowhere.
Let $$\cN_\Om(t)=\card(t\Omega\cap\bbZ^d),
$$
the number of integer lattice point inside the dilated domain $t\Omega$.
It is well known  that $\cN_\Om(t)$ is asymptotic
to $t^d\text{vol}(\Om)$
as $t\to \infty$. We denote by
$$
\Delta_\Om(t)=\frac{\cN_\Om  (t)-t^d\text{vol}(\Om)}{t^d\text{vol}(\Omega)}
\tag 1.1
$$
the relative error, or the discrepancy function.
It is conjectured that in dimensions $d\ge 5$  the relative error  is
$O(t^{-2})$ as $t\to\infty$. This conjecture is known to be true 
in the case of a ball centered at the origin, and for ellipsoids
in various degrees of generality
(see Landau \cite{20}, Walfisz \cite{31}, \cite{32}, Bentkus and G\"otze
\cite{2}). The error can be even smaller. For example,  
Jarn\'\i k \cite{14} established the bound $O(t^{-d/2+\eps})$ for the 
relative error, 
with  any $\eps>0$,
for almost all ellipsoids with axes parallel to the coordinate axes.
For general convex domains with non-vanishing curvature on the boundary,
W.  M\"uller \cite{22} proved that
$\Delta_\Om(t)= O( t^{-2+\lambda(d)+\eps}) $, where
$\lambda(d)=(d+4)/(d^2+d+2)$, if $d\ge 5$, $\lambda(4)=6/17$ and
$\lambda(3)=20/43$, improving on earlier results by Kr\"atzel and Nowak \cite{19}.
For planar domains, Huxley \cite{11} obtained this estimate with
$\lambda(2)=46/73$, which implies the relative error  
$O( t^{-100/73}(\log t)^{315/146})$.

In this paper we study the {\it mean square discrepancy}
of the lattice rest, the square function
$$
\cG_\Om(R)=\Big(\frac 1R \int_R^{2R}|\Delta_\Om(t)|^2 dt\Big)^{1/2}
\tag 1.2
$$
and related expressions.
For the ball $B_d$ in $\bbR^d$, centered at the origin,  bounds
and various asymptotics for mean square discrepancies 
have been obtained by
Walfisz \cite{32}  for $d\ge 4$,  Jarn\'\i k \cite{16}  for $d=3$  and Katai
\cite{17} for $d=2$.


In the more general situation where the boundary of $\Omega$ is smooth and is 
assumed to have everywhere non-vanishing Gaussian curvature, 
Nowak \cite{25} proved that
$\cG_\Om(R)=O(R^{-3/2})$ for planar domains.  This estimate is sharp by the
results of Bleher \cite{3} who  investigated the limit of 
$R^{3/2}\cG_\Om(R)$ as $R\to \infty$. 
The higher dimensional case was considered by
W. M\"uller \cite{21}, who proved a nearly sharp estimate for $d\ge 4$,
namely that
$\cG_\Om(R)\le C_\eps R^{-2+\eps}$ for any $\eps>0$. 
The case $d=3$ was left open.

The main  purpose of this note is to show that the known endpoint
bounds for the mean square
discrepancy in the case of the ball remain valid in the general case,
provided that $d\ge 4$. Moreover, we  prove a nearly sharp estimate
in dimension $d=3$, where we are off by a factor of $\sqrt {\log R}$.

\proclaim{Theorem 1.1} Let $\Omega$
be a convex domain in $\bbR^d$ containing the
origin in its interior, and assume that  $\Omega$ has  smooth boundary with
everywhere non-vanishing
Gaussian curvature.
Then there exists a constant $C(\Omega)$, such that for all $R\ge 2$,

$$
\cG_\Om(R) \le C(\Omega)\cases
R^{-2} \quad&\text{ if } d\ge 4
\\
R^{-2} \log R \quad&\text{ if } d=3
\\R^{-3/2}  \quad&\text{ if } d=2
\endcases .
\tag 1.3
$$
\endproclaim

As we noted above, the sharp estimate
$O(R^{-3/2})$ in the plane was already known for more general
{\it planar} domains with the non-vanishing curvature assumption.
In fact, it turns out that
this estimate
holds even if we
replace the mean square discrepancy over $[R, 2R]$ by
the mean square discrepancy over substantially smaller intervals
$[R, R+h]$.
A closely related result due to Huxley \cite{10} says that 
$(\int_R^{R+1}|\Delta_\Om(t)|^2 dt\Big)^{1/2} \leq  
C_\Om R^{-3/2}(\log R)^{1/2}$.

\proclaim{Theorem 1.2.1} Let $\Omega$
be a convex domain in $\bbR^2$ containing the
origin in its interior, and assume that  $\Omega$ has $C^\infty$ 
boundary with
non-vanishing  curvature.
Then there is a constant $C=C(\Omega)$ so that for all $R\ge 2$,
$$
\Big(\frac 1h \int_R^{R+h}|\Delta_\Om(t)|^2 dt\Big)^{1/2} \le
 C(\Om)
R^{-3/2} \quad \text { if }\quad \log R\le h\le R.
\tag 1.4
$$
\endproclaim


As an
immediate consequence of Theorem 1.2.1, the
mean square discrepancy over $[R,R+h]$ is dominated  
by $ C(\Om)  R^{-3/2}(\log R)^{1/2} h^{-1/2}$
if $ 0\le h\le \log R$. In particular,
the aforementioned result of Huxley follows if we set $h=1$.

We now consider more general domains in the plane. We say that a convex 
domain is of  type at most $\fm$  if its boundary has  order of contact at most $\fm$  with every
tangent line. Thus if $\fm=2$ we recover the case of
everywhere non-vanishing curvature considered above.
It is known that the analogue of Theorem 1.2.1 may fail if the order  
is greater than $2$ ({\it cf.} \cite{26}, \cite{5}, and \cite{23}). 
However for almost all rotations the estimate remains true for the rotated domain. More precisely  we have the following result. 

\proclaim{Theorem 1.2.2}
Let $\Omega$
be a convex domain in $\bbR^2$ containing the
origin in its interior, and assume that the boundary is smooth and
of finite  type at most $\fm$, in the sense that the order of contact of 
$\partial \Omega$ with every tangent line is at most $\fm$. 
For $A\in SO(2)$, denote by $A\Omega$ the rotated domain $ \{Ax:x\in\Omega\}$.
Then for almost all rotations $A$, the inequality  (1.4)
holds for $A\Om$,
with the  constant $C_{A\Om}$ depending on $A$.
More precisely, the following hold.

(i) The maximal function 
$$\cC_\Om(A)=
\sup_{R\ge 2} R^{3/2} \sup_{\log R\le h\le R}
\Big(\frac 1{h} \int_R^{R+h}|\Delta_{A\Om}(t)|^2 dt\Big)^{1/2}
\tag 1.5
$$
belongs to the weak type space  $L^{\frac {2\fm-2}{\fm-2},\infty}(SO(2))$. 

(ii) Let $\Ga$ be the set of all points $P\in \partial \Om$ where the 
curvature vanishes, and for $P\in \Gamma$ assume that the curvature vanishes 
of order
$m_P-2$ ($\le \fm-2$). Let $n_P$ be the outer unit normal at $P$ and $v_P$ a unit tangent vector at $P$.
Then $\cC_\Om(A)<\infty$, if $A$ satisfies, for some $\epsilon>0$, 
the Diophantine condition
$$
\max_{P\in \Gamma}\sup\{|k|^{\frac{m_P}{m_P-2}-\epsilon} 
|\inn{k}{A^*v_P}|:   \dist(k,\Bbb Rn_P)\le 1\} >0.
\tag 1.6
$$

In particular the set
$\{ A\in SO(2):\cC_\Om(A)=\infty\}$ is of 
Hausdorff dimension $\le \frac{\fm-2}{\fm-1}$.
\endproclaim

It is likely that one can weaken the Diophantine condition and thus the  
 estimate for the upper bound of the Hausdorff dimension is presumably 
not sharp.
The latter theorem is related to
 the results by Colin de Verdi\'ere \cite{5} and 
Tarnopolska-Weiss \cite{30}
who proved similar statements about the maximal function
$A\mapsto \sup_{t\ge 1} t^{4/3}\Delta_{A\Omega}(t)$.
Also, Nowak \cite{24} 
has obtained some improved van der Corput type  bounds
$| \Delta_{A\Omega}(t)|\le C_A t^{-4/3-\delta}$ for suitable 
$\delta=\delta(\Om)>0$, again
under appropriate  Diophantine conditions on the rotation.

We remark that in Theorem 1.2.1  the smoothness assumption
 can be relaxed considerably; 
moreover  a slightly weaker variant of Theorem 1.2.2 
holds without any assumption on the boundary of the convex domain. These issues are taken up in 
the  sequel \cite{13} to this paper.

{\it Notation:} 
Given two quantities $A$, $B$ we write $A\lesssim B$ 
if there is an 
absolute positive constant, 
depending only  on the specific domain $\Omega$, so that $A\le C B$.
We write $A\approx B$ if $A\lc B$ and $B\lc A$. 

\head{\bf 2. Preliminaries}\endhead

We denote by $\Omega^*$ the polar set of $\Omega$, 
$$\Omega^*=\{\xi: \inn{x}{\xi}\le 1 \text{ for all } x\in \Om\},
\tag 2.1$$
and let $\rho^*$ be the
Minkowski functional  associated to $\Omega^*$;
i.e. $\rho^*$ is homogeneous of degree $1$ and satisfies $\rho^*(\xi)=1$ if
$\xi\in \partial \Om^*$.
Then, if $P_+(\xi)$ is the unique point in
$\partial\Om$ at which $\xi$ is an outer normal to
$\partial\Om$, then
$$\rho^*(\xi)= \inn{P_+(\xi)}{\xi}. \tag 2.2.1
$$

Similary, if 
 $P_-(\xi)$ is the unique  point in
 $\partial\Om$ at which $\xi$ is an inner normal, then
$$\rho^*(-\xi)= -\inn{P_-(\xi)}{\xi}.
\tag 2.2.2
$$

If $t\mapsto x(t)$ is a regular $C^k$ parametrization of $\partial\Om$
near a point $P_0=x(t_0)$, 
and $t\mapsto n(t)$ denotes the outward unit normal vector, then
$t\mapsto x^*(t)=\inn {x(t)}{n(t)}^{-1} n(t)$ parametrizes
the boundary of $\Omega^*$, and $x^*$ is of class $C^{k-1}$.
If $\ka(P_0)$ denotes the Gaussian curvature at
$P_0$, and $\ka(P_0)\neq 0$ then the parametrization $t\mapsto x^*(t)$  
is regular near $P_0^*
=x^*(t_0)$ and the curvature $\kappa^*(P_0^*)$ of $\partial \Omega^*$
 at $P_0^*$
satisfies
$$|\kappa(P_0^*)\kappa(P_0)|=(|P_0| \cdot|P_0^*|)^{-d-1}.
\tag 2.3
$$
For these facts see e.g. Lemma 1 in \cite{21}.

We shall also need asymptotics for  the
indicator function of a convex domains. Suppose that
$\Omega$ is of finite line type
(in the sense that every tangent line has finite order of contact with $\partial \Omega$). Let $d\mu$ be a smooth density on the boundary of $\Omega$.
We define the Fourier transform by
$\widehat f(\xi)=\int f(y) \exp(\ic\inn{-y}{\xi}) dy$, and then a result
by Bruna, Nagel and Wainger
\cite{4} says that
$$\widehat {d\mu}(\xi)=
e^{-\ic\inn{P_+(\xi)}{\xi}}a_+(\xi)
+ e^{-\ic\inn{P_-(\xi)}{\xi}}a_-(\xi),
\tag 2.4.1
$$
where $a_{\pm}$ is smooth  and satisfies the symbol estimates
$$|\partial_\xi^{\alpha} a_\pm(\xi)|\le C_\alpha \gamma_\pm(\xi)
|\xi|^{-|\alpha|}, \quad |\xi|\ge 1
\tag 2.4.2
$$
for all multiindices $\alpha$, and $\gamma_\pm$ is defined as follows.
Let $H_{P}$ be the (affine) tangent plane to $\Omega$ at $P$.
Then  $\gamma_\pm (\xi)$ is the surface measure of the
cap
$$\gamma_\pm(\xi)=\sigma\big(\{y\in \partial \Omega:\dist(y,H_{P_\pm(\xi)})\le |\xi|^{-1}\}\big),
\tag 2.5
$$
where $\sigma$ denotes surface measure on $\partial \Om$.
By the divergence theorem, $\partial_{x_i}\chi_\Omega=-n_i d\sigma$, in
the  sense of distributions, where $n=(n_1,\dots, n_d)$ is the outward
 unit normal. Thus we get $\widehat {\chi_\Omega}(\xi)=-\ic\sum_{i=1}^d
(\xi_i/|\xi|^2)\widehat{n_i d\sigma}(\xi)$. If one 
combines  this with (2.2.1/2) and (2.4.1/2), one obtains
$$\widehat {\chi_\Omega}(\xi)=
e^{-\ic\rho^*(\xi)}b_+(\xi) +
e^{\ic\rho^*(-\xi)}b_-(\xi),
\tag 2.7
$$
where
$$|\partial_\xi^{\alpha} b_\pm(\xi)|\le C_\alpha
\gamma_\pm(\xi)
|\xi|^{-1-|\alpha|}, \quad |\xi|\ge 1.
\tag 2.8
$$

In the case of non-vanishing curvature
one has $\gamma_\pm(\xi)\lc|\xi|^{-(d-1)/2}$ but of course the above
statement, and more precise asymptotics, follow from the method of
stationary phase as in papers
by Hlawka \cite{8} (see also \S7 in \cite{9}).
More generally, for finite type domains one has 
$$\gamma_\pm(\xi)\lc |\kappa(x_\pm(\xi))|^{-1/2}|\xi|^{-(d-1)/2}. \tag 2.9
$$

This is proved in \cite{29}, and can also be deduced from the cap estimates (2.5)
using an argument in \cite{6}. However, it should be noted that these results are much easier in the two-dimensional case needed here. See \cite{27} and also \cite{1}.

%

\subheading{Definitions}
Let $\delta_0>0$ be fixed so that the ball $B_{2\delta_0}(0)$
with center $0$ and radius
$2\delta_0$ is contained in $\Omega$.
Let $\zeta$ be a smooth nonnegative radial
 cutoff function supported in the ball $B_\delta(0)$
 so that $\int\zeta(x) dx=1$. Let $\zeta_\eps(x)=\eps^{-d}\zeta(x/\eps)$.

We set $N(t)=\cN_\Om(t)$, 
$$E(t)=N(t)-t^d\volume (\Om),$$
 and

$$
\aligned
&N_\eps(t)=\sum_{k\in \bbZ^d} \chi_{t\Omega} * \zeta_{\eps} (k)\\&
E_\eps (t)=N_\eps (t)-t^d\volume(\Omega ).
\endaligned
\tag 2.10
$$
We also  denote by
$N^*_\eps(t)$ and $E^*_\eps(t)$ the corresponding expressions for the polar domain
$\Om^*$.


\subhead{\bf Three elementary Lemmas}  \endsubhead

\proclaim{Lemma 2.1} Suppose that $\Omega$ has $C^1$ boundary.
 Then there is a constant
$C=C(\Om)$ such that
for
$1\le R\le t\le 2R$, $0< \eps\le 1$,
$$
\alignat2
| E_\eps(t-\eps)|-Ct^{d-1}\eps&\le &|E(t)|&\le
|E_\eps(t+\eps)|+Ct^{d-1}\eps
\tag 2.11.1
\\
| E(t-\eps)|-Ct^{d-1}\eps&\le &|E_\eps(t)|&\le
|E(t+\eps)|+Ct^{d-1}\eps
\tag 2.11.2
\endalignat
$$

\endproclaim
\demo{\bf Proof}
By the properties of the the cutoff $\zeta_\eps$ we have
$$N_\eps (t-\eps )\le N(t)\le N_\eps(t+\eps ),$$
 and if we subtract $V(t)=t^d\volume(\Omega)$ throughout, we get
$$E_\eps (t-\eps)+[ V(t-\eps)-V(t)]\le E (t)\le
E_\eps (t+\eps)+[V(t+\eps)-V(t)].
$$
Clearly $|V(t\pm\eps)-V(t)|\lc t^{d-1}\eps$ and (2.11.1) follows. (2.11.2) follows as well if we apply (2.11.1)
with $t\pm\eps$ in place of $t$.\qed
\enddemo

\proclaim{Lemma 2.2} Suppose that $\mu\in[0,1]$ and that  the estimate
$$\sup_{t>0}t^{-(d-1-\mu)}|E(t)|\le C_1\tag 2.12$$ holds.
Then for $t\ge 1$ 
$$E_\eps(t)\lc \max\{t^{d-1-\mu}, t^{d-1}\eps\}.\tag 2.13
$$
Moreover there  a constant $C$ so that for 
$0<\eps\le h\le r$
$$
\Big|
\Big(\frac 1{h}\int_{r}^{r+h}|E(t)|^2 dt\Big)^{1/2}-
\Big(\frac 1{h}\int_{r}^{r+h}|E_\eps(t)|^2 dt\Big)^{1/2 } \Big|\le C
r^{d-1}
[h^{-1/2}\eps^{1/2} r^{-\mu}+
\eps].
$$


\endproclaim

\demo{\bf Proof} We first observe that (2.13) is immediate 
 by Lemma 2.1.
We integrate and obtain
$$
\align
\int_r^{r+h} |E(t)|^2 dt &\le
\int_r^{r+h+\eps} |E_\eps(t)|^2 dt+C h r^{2d-2}\eps^2
\\& \le
\int_r^{r+h} |E_\eps(t)|^2 dt+ C h r^{2d-2}\eps^2+C'\eps r^{2(d-1-\mu)}, 
\endalign
$$
which implies one of the desired inequalities, the other is obtained in the same way.\qed\enddemo

\proclaim{Lemma 2.3} Let $0<\eps<1$ and let for $\tau\ge 1$
$$\fS(\tau,\eps)=\card\{\ell\in\bbZ^d
:\tau-\eps\le\rho (\ell )\le \tau+\eps\}.
\tag 2.15$$
Then 
$$\fS(\tau,\eps)\le
C_1\tau^{d-1}\eps+C_2\Big(\int^{\tau+\eps/2}_{\tau-\eps/2} E_\veps(t)^2
N'_\veps (t)dt\Big)^{1/3}, \qquad \veps=\frac{4\eps}{\delta_0}.
\tag 2.16$$
\endproclaim

\demo{\bf Proof}
Let $t\in (\tau-\eps,\tau+\eps)$.
We use the elementary inequality
$$ \int\chi_{(t+h)\Om\setminus t\Om}(x-y)\veps^{-d}\zeta(\veps^{-1}y) dy \ge c_0 h/\veps
\quad\text{ if $h\ll\eps$, $x\in (\tau+\eps)\Om\setminus (\tau-\eps)\Om$.}
$$
This implies 
$$\align 
N_\veps(t+h)-N_{\veps}(t)&=
\sum_k 
\int\chi_{(t+h)\Om\setminus t\Om}(k-y)\veps^{-d}\zeta(\veps^{-1}y) dy 
\\&\ge  c_0 \frac{h}{\veps} \fS(\tau,\eps)
\endalign
$$ and thus
$$N_\veps'(t)\ge c_0 \fS(\tau,\eps) \eps^{-1},\quad |t-\tau|\le \eps, \quad
\veps=\frac{4\eps}{\delta_0}.
\tag 2.17
$$

We now turn to the proof of (2.16). We may assume that
$\fS(\tau,\eps)\ge C_1 \eps \tau^{d-1}$ with 
$C_1=d2^{d+1} c_0^{-1}\volume(\Om)$.
Then by (2.17), 
$$
\align
E_\veps'(t)&=N_\veps'(t)-  d\, t^{d-1}\volume(\Om)\\
&\ge N_\veps'(t)-  d\, (2\tau)^{d-1}\volume(\Om)\\
&\ge c_0 \fS(\tau,\eps) \eps^{-1} -
2^d C_1^{-1}d\,\eps^{-1}\volume(\Om)\fS(\tau,\eps)\\
&\ge c_0(2\eps)^{-1}\fS(\tau,\eps) .
\tag 2.18
\endalign
$$

Let $I_{\tau,\eps}=[\tau-\eps/2,\tau+\eps/2]$ and pick  $t_0\in I_{\tau,\eps}$ so that
$\min_{t\in I_{\tau,\eps}}|E_\veps(t)|=|E_\veps(t_0)|$; thus
$|E_\veps(t)|\ge|E_\veps(t)-E_\veps(t_0)|/2$ and
$|E_\veps(t)|\ge |\int_{t_0}^tE_{\veps}'(s)ds|/2 \ge c_0(4\veps)^{-1}|t-t_0|
\fS(\tau,\eps)$.
We use also (2.17) and obtain that
$$\align
\int^{\tau+\eps/2}_{\tau-\eps/2} E_\veps (t)^2 N'_\veps(t)dt
&\ge \int^{\tau+\eps/2}_{\tau-\eps/2}
(\tfrac{c_0 }{4\eps}\fS(\tau,\eps))^2|t-t_0|^2
\tfrac{c_0 }{\eps}\fS(\tau,\eps)\, dt\ge c [\fS(\tau,\eps)]^3
\endalign
$$
as asserted.\qed
\enddemo

\head{\bf 3. Proof of Theorem 1.1}\endhead
In this section we assume that $\Omega$ has a smooth boundary with
everywhere non-vanishing curvature. This implies that $\Omega^{*}$ is also 
smooth and has everywhere non-vanishing Gaussian curvature. See (2.3) above.   
We estimate the square-function 
$$G_\eps(R)=\Big(\frac 1R\int_R^{2R}|E_\eps(t)|^2 dt\Big)^{1/2}$$ for
$0<\eps\le 1/2$ and $R\ge 2$, and set
$$w_d(R)=\cases R^{2-d} &\quad\text{ if } d\ge 4\\
                (R\log R)^{-1} &\quad \text{ if } d=3\\
                R^{-1/2} &\quad \text{ if } d=2
\endcases
\tag 3.1
$$
and for $0<s\le 1/2$ let 
$$
A_d(s)=\sup_{s<\eps\le 1/2} \sup_{R\ge 2}
(1+\eps R)^{-d-1}
w_d(R) G_\eps(R).
\tag 3.2
$$
Analogously, we denote by
$A^*_d(s)$ the corresponding quantity associated to $\Om^*$.  It is not hard 
to see that $A_d(s)$ is finite 
for every $s$ since we have a trivial estimate
$A_d(s)\lc \sup_{R\ge 2} (1+s R)^{-d-1} R\lc s^{-1},$
and, similarly, 
$A_d^*(s)\lc s^{-1}$  for every
$s\le 1/2$. We shall see that $A_d(s)$ is bounded as $s\to 0$.  
Once 
this is established, the required bound for $\cG_\Omega$ follows from
$$\cG_\Omega(R)\lc R^{-d} \big(G_{1/R}(R)+ R^{d-2}\big)
, \tag 3.3$$
which is a consequence of Lemma 2.2.

The boundedness of $A_d(s) $  can be deduced from the following iterative 
procedure. 

\proclaim{Proposition 3.1}
There is a constant $C_\Om$ so that for $s\le 1/2$
$$A_d(s)^2 \le C_\Om\big(1+A_d^*(s)\big).
\tag 3.4$$
\endproclaim

Indeed, since $\Omega^{**} =\Omega$, (3.4) implies that   
$A_d^*(s)^2 \le C_{\Omega^*}\big(1+A_d(s)\big)$, so 
$$
A_d(s)^2\le C_\Om(1+\sqrt{C_{\Om^*}(1+A_d(s))})
$$
from which the boundedness of $A_d$ is immediate. 

\subheading{Proof of Proposition 3.1}
We estimate $G_\eps(R)$ assuming first that
$$R^{-1}\le \eps \le 1/2.$$
We apply  the Poisson summation formula $\sum_{k\in\bbZ^d} f(k)=
(2\pi)^d\sum_{k\in\bbZ^d}\widehat f(2\pi k)$ to $f=\chi_\Om(t\cdot)* \zeta_\eps$. This yields
$$E_\eps   (t)=\sum_{k\neq 0} (2\pi t)^d\widehat{\chi_\Om}
(2\pi tk)\widehat\zeta (2\pi\eps   k).
\tag 3.5.1
$$
We split
$E_\eps  (t)=\sum_{\pm} E^{\pm}_\eps  (t)$ by using (2.7/8); here
$$\aligned
E^{+}_\eps  (t)&=\sum_{k\neq 0} (2\pi t)^d b_+(2\pi tk)\exp(- 2\pi \ic\rho^*(k))\\
E^{-}_\eps  (t)&=\sum_{k\neq 0} (2\pi t)^d b_-(2\pi tk)
\exp( 2\pi \ic\rho^*(-k)).\endaligned 
\tag 3.5.2
$$

Now fix a nonnegative  $ \eta\in C^\infty(\Bbb R)$   so that  $\eta(t) =1$ for
$t\in [1,2]$ and $\eta$ is supported in $ (1/2, 3)$.
Then
$$ \align 
G_\eps(R)&\le G_\eps^+(R)+ G_\eps^-(R)
\\
&:=  \sum_\pm\Big(R^{-1}\int
|E^\pm_{\eps  } (t)|^2 \eta(R^{-1}t)dt\Big)^{1/2}.
\endalign
$$
We shall only consider estimates for $G_\eps^+(R)$ 
because the estimates for $G_\eps^-(R)$ are exactly analogous.
Multiplying out the squared expression we get
$$
G^+_\eps(R)^2=\sum\Sb k\neq 0\\ k'\neq 0\endSb
\widehat\zeta(2\pi\eps   k)\overline{\widehat\zeta(2\pi \eps   k' )}
R^{-1}\int e^{2\pi\ic t(\rho^*(k)-\rho^*(k'))}
q_{k,k'}(t)
dt 
\tag 3.6
$$ where
$$q_{k,k'}(t)=
b_+(2\pi tk)\overline{b_+(2\pi tk')} t^{2d}\eta(t/R) .
\tag 3.7
$$
Thus $q_{k,k'}$ is supported in $[R/2, 3R]$ and by (2.8) and $\gamma_\pm(\xi)=
O(|\xi|^{-(d-1)/2})$ we have the
symbol estimates $$\Big|\Big(\frac{d}{dt}\Big)^m q_{k,k'}(t)\Big|\le C_m
R^{d-1-m}|k|^{-(d+1)/2}
|k'|^{-(d+1)/2}.
\tag 3.8
$$
We now integrate by parts in $t$. We note that
 $|k| \approx\rho^*(k)$ and
$|\widehat\zeta(2\pi k/R)|\le C_N(1+|k/R|)^{-N}$ and
obtain the estimate
$$
G_\eps^+(R)^2\le C_{M,N} \sum_{k\neq 0}\sum_{k'\neq 0}
R^{d-1}
(1+R|\rho^*(k)-\rho^*(k' )|)^{-M}
(1+\eps|k|+\eps|k'|)^{-N}[\rho^*(k)\rho^*(k')]^{-\frac{d+1}{2}}.
$$

The terms with $|\rho^*(k)-\rho^*(k')|\ge R^{-1/2}$ give a contribution of 
$O(R^{d-1-M/2} \eps^{-2d})=
O(R^{3d-1-M/2})$ and we may choose $M=6d$.

Thus
$$\align
G_\eps^+(R)^2&\le C_1\sum_{-R^{1/2}\le n\le R^{1/2}}
\sum_{k\neq 0} (1+\eps\rho^*(k))^{-N}
\sum\Sb|\rho^*(\ell )-\rho^*(k)|\\ \in[\frac{n-1}R,\frac{n}{R}]\endSb
\frac{R^{d-1}}{(1+n)^M} [\rho^*(k)]^{-d-1} +C_2 R^{3d-1-M/2}
\\
&\le C_1' R^{d-1}\sum_{-R^{1/2}\le n\le R^{1/2}}(1+n)^{-M}
\sum_{k\neq 0}(1+\eps\rho^*(k))^{-N}
\frac{\fS^*(\rho^*(k),\tfrac {n+1}R)}
{\rho^*(k)^{d+1} }
+C_2' R^{3d-1-M/2};
\tag 3.9
\endalign
$$
here recall that $\fS^*(\tau, \eps)=\card\{\ell\in \Bbb Z^d:\tau-\eps\le \rho^*(\ell)\le \tau+
\eps\}$.
Now
$$
\align
&\sum_{k\neq 0} (1+\eps\rho^*(k))^{-N}
\frac{\fS^*(\rho^*(k), \tfrac {n+1}R)}
{\rho^*(k)^{d+1} }
\\&\lc
\sum_{l=0}^\infty 2^{-l}(1+\eps2^{l})^{-N}
\Big(\frac{1}{2^{ld}}\sum_{2^l\le \rho^*(k)<2^{l+1}}
[\fS^*(\rho^*(k),\tfrac {n+1}R)]^2\Big)^{1/2}\\
&\lc
\sum_{l=0}^\infty 2^{-l}(1+\eps2^{l})^{-N}
 [(n+1)I_{l}+II_{n,l}]
\endalign
$$
where
$$
\align I_{l}&= \Big(\frac{1}{2^{ld}}\sum_{2^l\le \rho^*(k)<2^{l+1}}
\rho^*(k)^{2d-2}R^{-2}\Big)^{1/2},
\\
II_{n,l}&=
\Big(\frac{1}{2^{ld}}
\sum_{2^l\le \rho^*(k)<2^{l+1}}
\int_{\rho^*(k)-(n+1)/2R}^{\rho^*(k)+(n+1)/2R}
\frac{{ {E_\veps^*}  (t)^2 N_\veps  ^*}'(t)}
{\fS^*(\rho^*(k),\tfrac {n+1}R)} dt\Big)^{1/2},
\endalign
$$
with $\veps=4\eps/\delta_0$;
here we used Lemma 2.3.
Observe that for $N$ large,
$$\align
\sum_{l=0}^\infty 2^{-l}(1+\eps 2^{l})^{-N}
I_l
\lc R^{-1}\sum_{l=0}^\infty 2^{l(d-2)}(1+\eps 2^{l})^{-N}
\lc \cases
R^{-1}\eps^{2-d}
&\quad \text{  if } d\ge 3
\\
R^{-1} \log (2+\eps^{-1})
&\quad \text{  if } d=2
\endcases
\tag 3.10
\endalign
$$
and thus, since we are assuming $\eps\le 1/R$,
$$
R^{d-1} \sum_{l=0}^\infty 2^{-l}(1+\eps 2^{l})^{-N} I_l\lc
R^{d-2}\max\{ \eps^{2-d} , \log (2+\eps^{-1})\}
\lc  w_d(R)^{-2}.
$$

We now estimate $II_{n,l}$ and set
 $J_{k,n}:=
[\rho^*(k)-(n+1)/2R,\rho^*(k)+(n+1)/2R]$. Observe that
$$
\align
\fS(\rho^*(k),\tfrac {n+1}R) &=\card \{\ell :\rho^*(k)-\tfrac{n+1}{R}\le \rho^*(\ell )\le
\rho^*(k)+\tfrac{n+1}{R}\}
\\&\ge \card\{ \ell : t-\tfrac{n+1}{2R}\le\rho^*(\ell )\le t+\tfrac{n+1}{2R}\}
\qquad \text{ if } |t-\rho^*(k)|\le \tfrac{n+1}{2R};
\endalign
$$
which is saying that
$\fS(\rho^*(k),\tfrac {n+1}R)\ge 
\fS(t,\tfrac {n+1}{2R})$ if $t\in J_{k,n}$. Thus
$$
\sum_{2^l <\rho^*(k)\le 2^{l +1}}
\frac{\chi_{J_{k,n}(t)}}{\fS(\rho^*(k),\tfrac {n+1}R)} \le
\frac{1}{\fS(t,\tfrac{n+1}{2R})}
\sum_{k}
\chi_{J_{k,n}}(t)=1.
$$
Therefore
$$
\align
II_{n,l}^2&=2^{-ld}\int
\Big[\sum_{2^l\le \rho^*(k)<2^{l+1}} \chi_{J_{k,n}}(t)
\Big]
\frac{{E_\veps^*}(t)^2 {N_\veps  ^*}'(t) }
{\fS^*(\rho^*(k),\tfrac {n+1}R)} dt
\\&\le
\frac{1}{2^{l d}}\int^{2^{l+2}}_{2^{l-1}}E^*_\veps  (t)^2
{N_\veps  ^*}'(t)dt\\
&=
\frac{1}{2^{ld}}
\int^{2^{l+2}}_{2^{l-1}}E^*_\veps   (t)^2 {E_\veps  ^*}' (t)dt
\,+\,\frac{1}{2^{l d}}\int^{2^{l +2}}_{2^{l-1}}
E^*_\veps   (t)^2\frac d{dt}(vol(t\Omega)) dt\\
&\le \frac{1}{2^{l d}}\Big(\frac{[E^*_\veps   (2^{l +2})]^3}3-
\frac{[E^*_\veps  (2^{l-1} )]^3}3\Big)+\frac{C}{2^l}\int_{2^{l-1}}^{2^{l +2}}
E^*_\veps   (t)^2dt\\
&\lc\Big( 2^{l(2d-6+\frac{6}{d+1})} +2^{l(2 d-3)}\veps^3+
\frac{1}{2^l}
\int^{2^{l +2}}_{2^{l-1}} E^*_\veps   (t)^2 dt\Big);
\endalign
$$
here we have used the estimate
$|E^*_{\veps  }(t)|\lc
2^{l(d-2+\tfrac 2{d+1})}+2^{l(d-1)}\veps$, $t\approx 2^l$,  which by Lemma 2.1 is a consequence of the classical 
 estimate $|E^*(t)|=O( t^{d-2+\frac{2}{d+1}})$, $d\ge 2$.
Thus 
$$
\align
&
R^{d-1} \sum_{l=0}^\infty 2^{-l}(1+\eps 2^{l})^{-N} II_{n,l}
\\
&\lc  R^{d-1}\Big(\sum_{l=0}^\infty
(1+\eps 2^{l})^{-N}\Big[  2^{l(d-4+\tfrac 3{d+1})}+ 
2^{l( d-5/2)}\veps^{3/2}+ 2^{-l}
\Big(\frac 1{2^l}\int_{2^{l-1}}^{2^{l+2}}|E_\veps^*(t)|^2 dt\Big)^{1/2}
\Big]\Big)
\\
&\lc  R^{d-1}\eps^{3-d} + R^{d-1}\sum_{l=0}^\infty 2^{-l} 
(1+\eps 2^{-l})^{-N} \sum_{i=-1}^{1} G_\veps(2^{l-i})
\\
&\lc
 R^{d-1}\eps^{3-d}
 + R^{d-1}\sum_{l=0}^\infty 2^{-l} (1+\eps 2^{-l})^{-N}
\frac{(1+\veps 2^l)^{d+1}}{w_d(2^l)} \sup_{r\ge 0} 
\big\{(1+\veps 2^r)^{-d-1}G_\veps^*(2^r) w_d(2^r)\big\}.
\tag 3.11
\endalign
$$
Now since
$R^{d-1}\eps^{3-d}\lc w_d(R)^{-2}$ for $\eps\ge R^{-1}$ we have 
$$
\align
&R^{d-1}\sum_{l=0}^\infty 2^{-l} \big(w_d(2^l)\big)^{-1}
 (1+\eps 2^{-l})^{-N}
(1+\veps 2^{-l})^{d+1} 
\\
&\lc R^{d-1}\sum_{l=0}^\infty 2^{-l} \big(w_d(2^l)\big)^{-1}
(1+\veps 2^{-l})^{-N+d+1}
\\
&\lc
  w_d(R)^{-2}(1+\veps R)^{d-2},
\tag 3.12
\endalign
$$
where the third inequality follows in a straightforward manner from the definition
of $w_d$. It is precisely at this point where one needs to distinguish 
the cases $d=2$, $d=3$ and $d\ge 4$.
Combining the previous estimates (3.10), (3.11) with (3.12)  we  obtain 
for $s\le 1$ and $\max\{s,R^{-1}\}\le \eps\le 1/2$
$$
\align
\big[(&1+\eps R)^{-d-1}w_d(R) G_\eps^+(R)\big]^2
\\
&\lc 1+ 
(1+\veps R)^{-2d-2}w_d(R)^2 R^{d-1}\sum_{l\ge 0}(1+\eps 2^l)^{-N} \sum_{|n|\le R^{1/2}}(1+n)^{-3}
\big((n+1)I_l+II_{n,l}\big)
\\
&\lc  1+ \sup_{r\ge 0} 
\big\{(1+\veps 2^r)^{-d-1}G_\veps^*(2^r) w_d(2^r)\big\}
\tag 3.13
\endalign
$$
for $\veps=4\eps/\delta_0$
The same estimate holds for $G_\eps^-(R)$ and thus for $G_\eps(R)$. 
Consequently, 
since $\eps\approx \veps$, we have 
$$
\big[(1+\eps R)^{-d-1}w_d(R) G_\eps(R)\big]^2\le C(1+A_d^*(s))
\qquad \text { if } R^{-1}\le \eps\le 1/2
\tag 3.14
$$

The required estimate  for $\eps\le 1/R$ follows from a small modification. Namely we can use 
Lemma 2.2 to see that
$$
\align G_\eps(R)&\le C_1\Big[
\Big(\frac 1R\int_R^{2R}|E(t)|^2 dt\Big)^{1/2} +R^{d-2}\Big]
\\
&\lc
C_2\Big[
\Big(\frac 1R\int_R^{2R}|E_{1/R}(t)|^2 dt\Big)^{1/2} +R^{d-2}\Big].
\endalign
$$
Thus 
$$
\align
&(1+\eps R)^{-2(d+1)}w_d(R)^2 G_\eps(R)^2
\\
&\qquad\lc
w_d(R)^2\big[ R^{2d-4}+
G_{1/R}(R)^2\big]
\\
&\qquad 
\le C(1+A_d^*(s))\qquad
\text{ if }\quad s\le \eps\le R^{-1}.
\tag 3.15
\endalign
$$
The desired estimate (3.4) follows from  (3.14), (3.15).\qed

\

\head{\bf 4. Localized square functions in the plane}\endhead

In this section we give the simple  proof of Theorem 1.2.1.
We assume that $\Omega$ is a convex domain in the plane, with smooth boundary, and that the curvature does not vanish at the boundary.

We may apply Lemma 2.2 with $\mu=0$, say,  and
we let $1\le h\le R$ and $\eps=R^{-1}$. Then
$$\Big(\frac 1h\int_R^{R+h}|E(t)|^2 dt\Big)^{1/2}
\lc\Big[
\Big(\frac 1h\int_R^{R+h}|E_{1/R}(t)|^2 dt\Big)^{1/2}+ (R/h)^{1/2}\Big].
\tag 4.1
$$
Let $\eta_0$ be a nonnegative $C^\infty$ function supported in $(-1/2,3/2)$
and which equals $1$ on $[0,1]$. Then
$$
\frac 1h\int_R^{R+h}|E_{1/R}(t)|^2 dt\lc\sum_{\pm}
\frac 1h\int |E^\pm_{1/R}(t)|^2\eta_0(\tfrac{t-R}{h}) dt
\tag 4.2
$$
with $E^\pm$ as in (3.5.2).
The  expressions on the right hand side are
 estimated by integration by parts, as in the previous section. We  square the series.
The cutoff  $\eta(t/R)$ is now replaced by
$\eta_0(\tfrac{t-R}{h})$ 
and this affects the argument since in the 
symbol estimates for 
the modification of $q_{k,k'}$ the estimate $R^{d-1-m}$ in (3.8) 
is now replaced 
by $R^{d-1} h^{-m}$.

As a result we obtain the estimate
$$\align
&\frac 1h\int |E^\pm_{1/R}(t)|^2\eta_0(\tfrac{t-R}{h}) dt
\\&\qquad \lc
R \sum_{k\neq 0}\sum_{k'\neq 0}  (1+h|\rho^*(k)-\rho^*(k')|)^{-M} 
(1+|k|/R+|k'|/R)^{-N}
|\rho^*(k)\rho^*(k')|^{-3/2}
\endalign
$$
and this term is estimated by a constant times
$$
\sum_{|n|\le R^{1/2}} \sum_{k \not=0}  
{(1+|n|)}^{-M} \rho^*(k)^{-3} {(1+\rho^{*}(k)/R)}^{-N} 
\fS^*(\rho^*(k)+
\tfrac nh,\tfrac 1h) +  R^{1-M/2},
\tag 4.3
$$
where, as before, 
$\fS^*(\tau,\eps)= 
\card\big(\{\ell\in \Bbb Z^2:\tau-\eps\le \rho^*(\ell)\le \tau+\eps\}\big).
$

Now by the classical estimate for the remainder term $E(t)$ with 
$t=\rho^*(k)+(n\pm 1)/h\approx \rho^*(k)$ we have
$$
\fS^*(\rho^*(k)+\tfrac nh,\tfrac 1h) \lc h^{-1}\rho^*(k)+\rho^*(k)^{2/3}.
\tag 4.4
$$ 

Putting the previous estimates together, we have
$$
\align
\frac 1h\int |E^\pm_{1/R}(t)|^2\eta_0(\tfrac{t-R}{h}) dt
&\lc
R\sum_{k\neq 0}(1+R^{-1}\rho^*(k))^{-N} \min\{ h^{-1}\rho^*(k)^{-2},
\rho^*(k)^{-7/3}\} + R^{1-M/2}
\\
&\lc R\big(1+h^{-1}\log R \big)
\endalign
$$
which is $O(R)$ if $h\gc \log R$. This finishes the proof of Theorem 1.2.1.\qed

\head{\bf 5. Estimates for finite type domains in the plane}\endhead
We shall give a proof of Theorem 1.2.2.
Let $\Omega$ be a convex finite type domain in $\Bbb R^2$ 
which contains the origin in its interior. 
We first give a version of the standard lattice rest estimate for the
polar set $\Omega^*$ which has a $C^1$ boundary.

\proclaim{Lemma 5.1} We have the following estimate for the Fourier transform 
of  the characteristic function of $\Omega^*$,
$$
\big|\widehat {\chi_{\Om^*} }(\xi)\big|\le C(1+|\xi|)^{-3/2}.
\tag 5.1
$$
\endproclaim

Taken Lemma 5.1 for granted we obtain as a consequence
\proclaim{Corollary 5.2}
 Let $\Omega$ be a convex set in $\Bbb R^2$, containing the origin in 
its interior and suppose that $\Omega$ has smooth finite type boundary. Let $\Omega^*$ be the polar set.
Then
$$
\cN_{\Om^*} (t) = t^2 \text{area}(\Omega^*) +O(t^{2/3})
\tag 5.2
$$
as $t\to\infty$.
\endproclaim
\demo{Proof}
This follows from Lemma 5.1 using the standard argument (see e.g. \cite{8}, or
\S7 of  \cite{9}).
\enddemo

The Corollary can be improved by using more sophisticated techniques
which however are not needed here.

Before proving Lemma 5.1 we recall some terminology:
We denote by $\Ga$  the set of
all points in $\partial\Omega$ 
at which the curvature vanishes; these points are separated and thus $\Ga$ is finite. For every $P\in \Ga$ 
let $m_P$ be the type  at $P$
(i.e. 
 the curvature vanishes of order $m_P-2$ at $P$). 
For every $P\in \partial\Om$ there is a unique $P^*\in \partial \Omega^*$ 
so that $\inn{P}{P^*}=1$ and we define  $\Gamma^*=\{P^*: P\in \Gamma\}$.

\demo{\bf Proof of Lemma 5.1}

The boundary $\partial \Omega^*$  is smooth away from $\Gamma^*$ and it is $C^1$ everywhere.
Thus surface measure $d\sigma$ is  well defined and by 
an application of the divergence theorem as 
in \S2 
estimate  (5.1) follows provided  we can show that
$$
|\widehat {\chi d\sigma}(\xi)|\lc (1+|\xi|)^{-1/2}
\tag 5.3
$$
for $\chi\in C^\infty_0$.

To see this we introduce a partition of unity $\chi d\sigma=\sum_\nu
\chi_\nu d\sigma$ 
where  each $P^*\in \Gamma^* $  lies in exactly one of the supports of the 
functions $\chi_\nu$.
Clearly it suffices to prove the estimate $\widehat {d\sigma_\nu}(\xi)=
O(|\xi|^{-1/2})$  for each $\sigma_\nu:=\chi_\nu d\sigma$.


Fix $\nu$ and $P\in \Gamma$. If $P^*\notin \supp d\sigma_\nu$ then 
$\widehat {d\sigma_\nu}(\xi)=
O(|\xi|^{-1/2})$  by the standard stationary phase argument.
Thus suppose $P\in \Gamma\cap \supp d\sigma_\nu$.
 By a rotation we may assume that $n_P=(0,1)$ and by an additional 
translation we may also assume that $P$ lies on the $x_2$-axis. Let $m=m_P$ be the type at 
$P$. Near $P$ the boundary of $\Omega$ is parametrized by $(t, f(t))$ where 
$$f(t)= a_0-a_m\frac{t^m}{m}+ t^{m+1} g_1(t)$$
with $a_0>0$, $a_m>0$.
Thus  a parametrization of $\partial \Omega^*$ 
 near $P^*=(a_0^{-1},1)$ is given by
$$
t\mapsto 
\frac{n(t)}{\inn{x(t)}{n(t)}}= 
\frac 1{f(t)-tf'(t)}
\frac{(-f'(t),1)}{\sqrt{1+f'(t)^2}};
$$
however this parametrization is not regular.
Denote by
$\omega(t)$ the first coordinate of 
${\inn{x(t)}{n(t)}}^{-1}n(t)$. Then it is easy to see that
$$
\omega(t)= (a_m/a_0)
t^{m-1}(1+ t g_2(t))= (c_0 s(t))^{m-1}
$$
where $c_0= (a_m/a_0)^{1/(m-1)}$ and $s(t)=t+O(t^2)$.
Moreover
$$ 
(f(t)-tf'(t))^{-1}(\sqrt{1+f'(t)^2})^{-1} =a_0^{-1}\big(
1-\tfrac{m-1}m\tfrac {a_m}{a_0} t^m +t^{m+1}g_2(t)\big).
$$
Thus setting $\tau= (c_0 s(t))^{m-1}$ we see after a short computation that
near $P^*$ the boundary is parametrized by
$\tau\mapsto (\tau, h(\tau))$ with
$$h(\tau)= a_0^{-1} \big(1-c_1 \tau^{m/(m-1)}+
\tau^{\frac{m+1}{m-1}} g_3(\tau^{\frac 1{m-1}})\big)
$$
where $c_1=(m-1)m^{-1}(a_m/a_0)c_0^{-m}=
(m-1)m^{-1}(a_m/a_0)^{-1/(m-1)}$ and $g_3$ is smooth.
Thus we have to show that
$$
J(\xi)=\int e^{-i(\xi_1 \tau+\xi_2 h(\tau))} \eta(\tau) 
d\tau =O( |\xi|^{-1/2})
\tag 5.4
$$
as $|\xi|\to \infty$; here we may 
assume that the support of $\eta_\nu$ is contained in a 
small interval $(-\delta,\delta)$.

It suffices to estimate the analogous integral extended over the set 
$\{\tau:|\tau|\ge|\xi|^{-1/2}$ since the error is $O(|\xi|^{-1/2})$.
Observe that for small $\tau$ we have $|h'(\tau)|\ll 1$ and $|h''(\tau)|\ge c
\tau^{-(m-2)/(m-1) }\gc 1$. 
Thus by van der Corput's lemma (\cite{28}, ch. VIII.1)
we obtain for large $|\xi|$ the estimate
$|J(\xi)|\lc |\xi_1|^{-1}$ if $|\xi_1|\ge |\xi_2|$ (using first derivatives of the phase function) and the estimate
$|J(\xi)|\lc |\xi_2|^{-1/2}$ if $|\xi_2|\ge |\xi_1|$ (using second 
derivatives).
This implies (5.4) and thus (5.3).\qed
\enddemo

\subheading{Proof of Theorem 1.2.2}
We shall decompose the Fourier transform of $\chi_\Om$ as in
 \cite{27}, following rather closely \cite{13}. 
Using the divergence theorem as above, we see that 
$$
\widehat {\chi_\Omega}(\xi)=
\ic \sum_{i=1}^d 
\frac{\xi_i}{|\xi|^{2} }  \int_\Si n_i(y)e^{-\ic\inn{y}{\xi}}d\sigma(y)
\tag 5.5
$$
where $n_i$ denotes the $i^{\text {th}}$ component of $n_P$. \

For every $P \in \Gamma$ we choose a 
narrow  conic symmetric  neighborhood $V_P$  of the normals $\{\pm n_P\}$, a small
neighborhood  $U_P$ of $P$  in $\Sigma$ and a $C^\infty_0$ 
function $\chi_P$ whose restriction to $\Si$ vanishes off $\cU$ 
and so that $\chi_P$ equals one in a neighborhood of $P$. 
We may arrange these neighborhoods so that the sets 
$\overline V_P\cap\{\xi:|\xi|\ge 1\}$,
$P\in \Gamma$ are pairwise disjoint and that  the normals to all points 
in a neighborhood of $\overline U_P$ are contained 
in $V_P$, so that  the  $\overline U_P$'s  are disjoint too.

Define 
$$F_{i,P}(\xi)=\int_\Si \chi_P(y)n_i(y)e^{-\ic\inn{y}{\xi}}d\sigma(y)
$$
Let $v_P$ a unit tangent vector to $\partial \Om$ at $P$.
Then if 
the cones  $V_P$ are chosen sufficiently narrow, we have

$$\sum_{i=1}^d \frac{\xi_i}{|\xi|^2} F_{i,P}(\xi) =
e^{-\ic\rho^*(\xi)}b_+(\xi)
+ e^{\ic\rho^*(-\xi)}b_-(\xi)
\tag 5.5
$$
where 
$$|\partial_\xi^{\alpha} b_\pm(\xi)|\le C_\alpha \cases |\xi|^{-1-|\alpha|}
\min\{|\xi|^{-\frac{1}{m_P}},\,  
\xi^{-\frac 12} 
\Theta_P(\xi)
\}
\quad &\text{ if }\xi\in V_P
\\
C_N|\xi|^{-N}&\text{ if } \xi \notin V_P;
\endcases 
\tag 5.6
$$
 with
$$\Theta_P(\xi) =\Big|\frac{\inn{v_P}{ \xi}}{\inn{n_P}{\xi}}
 \Big|^{-\frac{m_P-2}{2(m_P-1)}}.\tag 5.7
$$
This follows from (2.8) (with $\alpha=0$) and (2.9) by a 
straightforward computation.
Moreover
$$ 
\sum_{i=1}^d \frac{\xi_i}{|\xi|^2} \big(F_i(\xi)-\sum_{P\in \Gamma}F_{i,P}(\xi)\big)
 =
e^{-\ic\rho^*(\xi)}c_+(\xi)
+ e^{\ic\rho^*(-\xi)}c_-(\xi)
\tag 5.8
$$
where
$$
|\partial_\xi^{\alpha} c_\pm(\xi)|\le C_\alpha 
 |\xi|^{-3/2-|\alpha|}. \tag 5.9
$$

The estimate for $\xi\in V_P$ follows from Proposition 1.2, and the estimate for $\xi\notin V_P$ follows by a simple integration by parts;
namely if $t\mapsto \gamma(t)$ parametrizes $\Si$ near $P$ 
then $|\inn{\gamma'(t)}{\xi}|\approx
|\xi|$ for $\gamma(t)\in U_P$ and $\xi\notin V_P$.

Moreover by the usual stationary phase or van der Corput  estimate we have 
$$
|F_i(\xi)-\sum_{P\in\Gamma} F_{i,P}(\xi)|\lc(1+|\xi|)^{-1/2}\tag 5.10
$$
here we used the definition of $\Gamma$ and the fact that  $\chi_P$  is equal
to $1$ near $P$.


Let $E_{1/R, A}(t)$ be the remainder term (2.10) with $\eps=1/R$, 
with $\Om$ replaced by  the
 rotated domain $A\Om$; that is
$$
\align
E_{1/R, A}(t)&=\sum_{k\in \bbZ^2} \chi_{t\Omega} * \zeta_{1/R} (A^{-1}k)-t^2
\text{area}(\Omega )
\\
&=\sum_{k\neq 0} t^2 \widehat\zeta(2\pi R^{-1} Ak) 
\sum_{i=1}^d \frac{2\pi t\inn{Ak}{e_i}}{|2\pi tAk|^2}
F_i(2\pi tAk)
\tag  5.11
\endalign
$$

For $P\in \Gamma$, $A\in SO(2)$ let
$$
\align
\cZ^P_I(A)&=\{k\in \Bbb Z^d: Ak\in V_P, k\neq 0, \dist (Ak, \bbR n_P)<1\}
\\
\cZ^P_{II}(A)&=\{k\in \Bbb Z^d: Ak\in V_P, k\neq 0, \dist (Ak, \bbR n_P)\ge 1\}
\endalign
$$
and let
$$
\cZ_{III}(A)=\{ k\in\bbZ^d: k\neq 0, k\notin \cup_{P\in\Gamma}V_P\}.
$$

We may use  estimate (4.1) which does not depend on any curvature assumptions and  see that
it suffices to estimate 
the square function
$(h^{-1} \int |E_{1/R,A}(t)|^2  \eta_0(\tfrac{t-R}{h}) dt)^{1/2}$ ({\it cf.} (4.2)). 
We decompose for $R\le t\le 2R$

$$\align E_{1/R, A}(t)=&\Big(\sum_{P\in \Ga}\sum_{k\in \cZ_I^P(A)}+
\sum_{P\in \Ga}\sum_{k\in \cZ_{II}^P(A)}
t^2 \widehat\zeta(2\pi R^{-1} Ak) \sum_{i=1}^d \frac{2\pi t\inn{Ak}{e_i}}{|2\pi tAk|^2}
F_{i,P}(2\pi tAk)
\\
&\quad
+
\sum_{k\in \cZ_{III}(A)}
t^2 \widehat\zeta(2\pi R^{-1} Ak) \sum_{i=1}^d \frac{2\pi t\inn{Ak}{e_i}}{|2\pi tAk|^2}
\big(F_i(2\pi tAk)-\sum_{P\in \Ga}
F_{i,P}(2\pi tAk)\big)
\\&\quad+
\sum_P\sum_{k\notin V_P} 
t^2 \widehat\zeta(2\pi R^{-1} Ak) \sum_{i=1}^d \frac{2\pi t\inn{Ak}{e_i}}{|2\pi tAk|^2}
F_{i,P}(2\pi tAk)
\\
=&
 \sum_{\pm}\Big(\sum_{P\in\Ga} I_P^\pm(t)+\sum_{P\in \Gamma}II_P^\pm(t)+III^\pm(t)\Big)+IV(t)
\endalign
$$
where $$|IV(t)|=O(t^{-N})
\tag 5.12$$ and 
$$
\align
I_P^+(t,A) &= \sum_{k\in \cZ^P_{I}(A)} \widehat\zeta(2\pi R^{-1} Ak) b_+(2\pi t Ak) e^{-2\pi i t \rho^*(Ak)}
\\
II_P^+(t,A) &= \sum_{k\in \cZ^P_{II}(A)} \widehat\zeta(2\pi R^{-1} Ak)  b_+(2\pi t Ak) e^{-2\pi i t \rho^*(Ak)}
\\
III^+(t,A) &= \sum_{k\in \cZ^P_{III}(A)}\widehat\zeta(2\pi R^{-1} Ak)  c_+(2\pi t Ak) e^{-2\pi i t \rho^*(Ak)},
\endalign
$$
and the expressions $I^-_P$, $II^-_P$ and $III^-_P$ are defined by replacing 
$b_+$ by $b_-$, $c_+$ by $c_-$, and 
$e^{-2\pi i t \rho^*(Ak)}$ by
$e^{2\pi i t \rho^*(-Ak)}$.

The argument in the previous section applies to the square functions associated to $III^\pm(t,A)$ and we obtain the bound
$$
\frac{1}{h}\int|III^\pm(t,A)|^2 
\eta_0(\tfrac{t-R}{h})
dt\lc R(1+h^{-1}\log R),
\tag 5.13
$$
uniformly in $A$.

A small variation of this argument also applies to the square function associated to $II_P^\pm(t,A)$. Namely, arguing as in \S3 and using (5.6/7) we see that

$$
\align
&\frac{1}{h}\int|II_P^+(t,A)|^2 \eta_0(\tfrac{t-R}{h})dt\\
\\&\lc 2
\sum\Sb k\in \cZ^P_{II}(A)\\ k'\in \cZ^P_{II}(A)\\ \Theta_P(Ak)\ge
 \Theta_P(Ak')\endSb
R(1+h|\rho^*(Ak)-\rho^*(Ak')|)^{-N}
\big(1+\frac{|k|+|k'|}R\big)^{-N} 
\frac{\Theta_P(Ak)\Theta_P(Ak')}{\rho^*(Ak)^{3/2}
\rho^*(Ak')^{3/2}}
\\
&\lc R\sum_{k\in \cZ^P_{II}(A)}
\frac{\Theta_P(Ak)^2}{\rho^*(Ak)^{3}}
\sum\Sb n\in \bbZ\\|n|\le 2^{k-5}\endSb(1+n)^{-N}
(1+|k|/R)^{-N}\fS_A^*(\rho_*(Ak)+\tfrac nh,\tfrac 1h)+R^{4-N}
\endalign
$$ 
where now
$\fS_A^*(\tau,\eps)= 
\card\big(\{\ell\in \Bbb Z^2:\tau-\eps\le \rho^*(A\ell)\le \tau+\eps\}\big).
$

Observe that  $\dist (Ak, \Bbb R n_P)\ge 1$ and 
 $\dist (A\xi, Ak)\le 1/2$  implies that
$\Theta_P(Ak)\approx \Theta_P(A\xi)$. 
Thus we can use the argument in 
\S3 and Lemma  5.1 
and estimate
$$\align
&\frac{1}{h}\int_R^{R+h}|II_P^\pm(t,A)|^2 \eta_0(\tfrac{t-R}{h}) 
dt\\ &\qquad\lc R
\int_{V_P}(1+|\rho^*(\xi)|)^{-3}\Theta_P^2(\xi)
(1+R^{-1}\rho^*(\xi))^{-N} \min\{ h^{-1}\rho^*(\xi)^{-2},
\rho^*(k)^{-7/3}\} d\xi + R^{1-M/2}.
\endalign
$$
Since $\Theta_P^2$ is homogeneous of degree $0$ and 
integrable over the sphere $\{\rho^*(\eta)=1\}$
it is easy to see that the former expression is bounded by
$R(1+h^{-1}\log R),$
thus
$$
\frac{1}{h}\int_R^{R+h}|II_P^+(t,A)|^2 dt\lc R(1+h^{-1}\log R),
\tag 5.14
$$
for $|h|\le R$, uniformly in $A$. The same estimate holds
true with $II_P^+$ replaced by $II_P^-$- the proof only requires changes in
the notation.

In order to estimate the square function involving
$I_P^+ $, we let $S_P(A)$ be the set of all  $k\in \Bbb Z^2\setminus\{0\}$
with $\dist(k, \Bbb R A^*n_P)<1$, and define 
$$
M_{P,\eps}(A)=  \sup
\big\{|k|^{-1+\eps}\Theta_P(k): k\in S_P(A)\big\} .
$$ 

Then
$$\align
&\frac{1}{h}\int|I^+_P(t)|^2\eta_0(\tfrac{t-R}{h}) dt
\\
&\lc \sum_{k\in S_P(A)}
\sum_{k'\in S_P(A)}
R(1+h|\rho^*(Ak)-\rho^*(A k')|)^{-N}(1+\tfrac{|k|+|k'|}R)^{-N}
\frac{\Theta_P(Ak)\Theta_P(Ak')}{|k|^{3/2}|k'|^{3/2}}
\\&\lc M_{P,\eps}(A)^2 
\sum_{k\in S_P(A)}
R(1+h|\rho^*(Ak)-\rho^*(Ak')|)^{-N}(1+\tfrac{|k|+|k'|}R)^{-N}|k|^{-\eps-1/2}
|k'|^{-\eps-1/2}
\\&\lc M_{P,\eps}(A)^2 
\sum_{k\in S_P(A)}
R(1+|k|/R)^{-N}|k|^{-1-2\eps},
\endalign
$$
and thus 
$$\frac 1h\int|I^+_P(t)|^2\eta_0(\tfrac{t-R}{h}) dt\le C_\eps M_{P,\eps}(A)^2
 R. \tag 5.15
$$
Again the same estimate remains true for 
$I^-_P(t)$.

For each $k\neq 0$ the function  $A\mapsto  \Theta_P(Ak)$ belongs 
to  the space $L^{(2m_P-2)/(m_P-2),\infty}$. For $\alpha>0$ the set
$\{A\in SO(2): M_{P,\eps}(A)>\alpha\}$ is the union of the sets
$E_k(\alpha)=\{A: \Theta_P(Ak)> |k|^{1-\eps}\alpha\}$, $k\in\Bbb Z^2\setminus \{0\}$
and the measure of $E_k(\alpha)$ is
$\lc (k^{1-\eps}\alpha)^{-(2m_P-2)/(m_p-2)}$. Since 
$(2m_P-2)/(m_p-2)>2$ we may sum over all  $k\in \Bbb Z^2\setminus\{0\}$ and we see that 
$M_{P,\eps}\in  L^{(2m_P-2)/(m_P-2),\infty}(SO(2))$ provided that $\eps\le 1/2$.
Combining the estimates (5.12-5.15) this proves
that $\cC_\Om\in L^{(2m_P-2)/(m_P-2),\infty}(SO(2))$.

The Diophantine condition (1.6) for some $\veps>0$  is equivalent with the condition
$M_{P,\eps}(A)<\infty$, for some $\eps>0$. Fix $P$.
The estimates (5.12-15) show that  $\cC_\Omega(A)=\infty$ also implies 
$M_{P,\eps}(A)=\infty$ for at least one $P\in \Gamma$. Thus we  can complete  the proof if
for any sufficiently small $\eps>0$ we demonstrate that 
the set $\{A\in SO(2): M_{P,\eps}(A)=\infty\}$ 
has Hausdorff dimension  $\le (m_P-2)(m_P-1)^{-1}(1-\eps)^{-1}$.

 Set $\beta=(m_P-2)/(2m_P-2)$, thus $\beta<1/2$.
Now $M_{P,\eps}(A)=\infty$ implies that there are infinitely many 
$k\in S_P(A)$  so that $|k|^{\eps-1}|\inn{k/|k|}{v_P}|^{-\beta}\ge 1$. 
If $A^*v_P=(\alpha_1,\alpha_2)$ this means 
$|k_1\alpha_1+k_2\alpha_2|\le |k|^{(\beta-1+\eps)/\beta}$.
Now $|\alpha_1|\ge |\alpha_2|$ implies $|k_1|\lc  |k_2|$
and
$|\alpha_2|\ge |\alpha_1|$ implies $|k_2|\lc  |k_1|$ (as $k\in S_P(A)$).
Thus if  $|\alpha_1|\ge |\alpha_2|$ 
 the condition $M_{P,\eps}(A)=\infty$ implies 
that  for infinitely many $k$ with $|k_2|\approx |k|$ we have that
 $$
|k_1/k_2-\alpha_2/\alpha_1|\le C |k_2|^{(\eps-1)/\beta}
\qquad\text{ or }\qquad
|k_1/k_2+\alpha_2/\alpha_1|\le C |k_2|^{(\eps-1)/\beta}.
\tag 5.16.1
$$
Likewise, 
 if   $|\alpha_2|\ge |\alpha_1|$ 
 and  $M_{P,\eps}(A)=\infty$ then 
$$
|k_2/k_1-\alpha_1/\alpha_2|\le C |k_1|^{(\eps-1)/\beta}
\qquad\text{ or }\qquad
|k_2/k_1+\alpha_1/\alpha_2|\le C |k_1|^{(\eps-1)/\beta}
\tag 5.16.2
$$
  for infinitely many $k$ with  $|k_1|\approx |k|$.

Let $P_\theta$ denote the set of all   $x\in [-1,1]$ 
for which there exists infinitely many rationals $p/q$ such that
$|x-p/q|\le q^{-2-\theta}$. By a Theorem of
Jarn\'\i k \cite{15} (see also \cite{18})  the Hausdorff dimension of $P_\theta$
is equal to $2/(2+\theta)$ (and we need only the easy upper bound).
Now choose in  (5.16.1/2) a small $\eps>0$
(in particular so that 
$\beta< (1-\eps)/2$)
 and we apply the last statement with
 $\theta=(1-\eps)\beta^{-1}-2$ and then
$2/(2+\theta)=2\beta(1-\eps)^{-1}=(m_P-2)(m_P-1)^{-1}(1-\eps)^{-1}$.

Consequently,  with $\fm$ being the maximal type,
the Hausdorff dimension of the  set $\{A\in SO(2):\cC_\Omega(A)=\infty\}$ 
does not exceed $(\fm-2)/(\fm-1)$.
\qed

\enddocument